\documentclass{birkmult}

\usepackage{amssymb}
\usepackage{latexsym}
\usepackage{amsmath}
\usepackage{xypic}
\usepackage{amsbsy}

 \newtheorem{thm}{Theorem}[section]
 \newtheorem{cor}[thm]{Corollary}
 \newtheorem{lem}[thm]{Lemma}
 \newtheorem{prop}[thm]{Proposition}
 \theoremstyle{definition}
 \newtheorem{defn}[thm]{Definition}
 \theoremstyle{remark}
 \newtheorem{rem}[thm]{Remark}
 
 \numberwithin{equation}{section}

\def\qed{\hfill$\Box$}

\def\CH{{\mathcal {H}}}
\def\CK{{\mathcal {K}}}
\def\CP{{\mathcal {P}}}
\def\CD{{\mathcal {D}}}
\def\CO{{\mathcal {O}}}

\newcommand{\eins}{\mathbf{1}}

\def\ua{{\underline{a}}}

\def\ut{{\underline{t}}}
\def\us{{\underline{s}}}
\def\uT{{\underline{T}}}

\def\uV{{\underline{V}}}

\def\uA{{\underline{A}}}
\def\uB{{\underline{B}}}

\def\N{{\mathbb {N}}}
\def\C{{\mathbb {C}}}

\newcommand{\aso}{\stackrel{\circ}{a}^*}
\newcommand{\Ao}{\AA}
\newcommand{\Aso}{{\AA}^*}
\newcommand{\Ho}{\stackrel{\circ}{\CH}}

\newcommand{\uAo}{\stackrel{\circ}{\uA}}
\newcommand{\uAso}{\stackrel{\circ}{\underline{A}^*}}
\newcommand{\Do}{\stackrel{\circ}{D}}
\newcommand{\Dso}{\stackrel{\circ}{D}_*}
\newcommand{\CDo}{\stackrel{\circ}{\CD}}
\newcommand{\CDso}{\stackrel{\circ}{\CD}_*}
\newcommand{\Po}{\stackrel{\circ}{\CP}}
\newcommand{\Co}{\stackrel{\circ}{C}}

\newcommand{\oomega}{\overline{\omega}}
\newcommand{\uomega}{\underline{\omega}}

\newcommand{\Od}{\mathcal{O}_d}
\newcommand{\OH}{\Omega_{\CH}}
\newcommand{\OK}{\Omega_{\CK}}
\newcommand{\OP}{\Omega_{\CP}}

\begin{document}

%
%
%
%
%
%
%
%
%
\title[Characteristic Functions for Ergodic Tuples]
 {Characteristic Functions for Ergodic Tuples}
\author{Santanu Dey}

\address
{
Institut f\"ur Mathematik und Informatik,\\
Ernst-Moritz-Arndt-Universit\"at,\\
Friedrich-Ludwig-Jahn-Str. 15a, \\
17487 Greifswald, Germany.
}

\email{dey@uni-greifswald.de}


\author{Rolf Gohm}
\address
{
Department of Mathematics, \\
University of Reading, \\
Whiteknights, P.O.Box 220, \\
Reading, RG6 6AX, England
}
\email{r.gohm@reading.ac.uk}



\date{June 17, 2005}

\begin{abstract}
Motivated by a result on weak Markov dilations, we define a notion
of characteristic function for ergodic and coisometric row contractions
with a one-dimensional invariant subspace for the adjoints. This extends a
definition given by G.\,Popescu. We prove that our characteristic function
is a complete unitary invariant for such tuples and show how it can be
computed.
\end{abstract}

\maketitle

\setcounter{section}{-1}
\section{Introduction}

If $Z = \sum^d_{i=1} A_i \cdot A^*_i$ is a normal, unital, ergodic,
completely positive map on $B(\CH)$, the bounded linear operators
on a complex separable Hilbert space, and if there is a (necessarily
unique) invariant vector state for $Z$, then we also say that
$\uA = (A_1,\ldots,A_d)$ is a coisometric, ergodic row contraction
with a one-dimensional invariant subspace for the adjoints. Precise
definitions are given below. This is the main setting to be
investigated in this paper.

In Section 1 we give a concise review of a result  on the dilations
of $Z$ obtained by R.\,Gohm in \cite{Go04} in a chapter called `Cocycles
and Coboundaries'. There exists a conjugacy between a homomorphic dilation
of $Z$ and a tensor shift, and we emphasize an explicit infinite product formula
that can be obtained for the intertwining unitary.
\cite{Go04} may also be consulted for connections of this topic to a
scattering theory for noncommutative Markov chains by B.\,K\"ummerer
and H.\,Maassen (cf. \cite{KM00}) and more general for the relevance of this
setting in applications.

In this work we are concerned with its relevance in operator theory
and correspondingly in Section 2 we shift our attention to the row
contraction $\uA = (A_1,\ldots,A_d)$. Our starting point has been the
observation that the intertwining unitary mentioned above has many
similarities with the notion of characteristic function occurring in the
theory of functional models of contractions, as initiated by
B.\,Sz.-Nagy and C.\,Foias (cf. \cite{NF70,FF90}). In fact, the center of
our work is the commuting diagram 3.3 in Section 3, which connects the
results in \cite{Go04} mentioned above with the theory of minimal isometric
dilations of row contractions by G.\,Popescu (cf. \cite{Po89a}) and
shows that the intertwining unitary determines a multi-analytic inner function, in the sense
introduced by G.\,Popescu in \cite{Po89c,Po95}. We call this inner function
the {\it extended characteristic function} of the tuple $\uA$, see Definition 3.3.

Section 4 is concerned with an explicit computation of this inner function.
In Section 5 we show that it is an extension of the characteristic function of the
$*$-stable part $\uAo$ of $\uA$, the latter in the sense of Popescu's
generalization of the Sz.-Nagy-Foias theory to row contractions (cf. \cite{Po89b}).
This explains why we call our inner function an {\it extended} characteristic function.
The row contraction $\uA$ is a one-dimensional extension of the $*$-stable row contraction
$\uAo$, and in our analysis we separate the new part of the characteristic function from
the part already given by Popescu.

G.\,Popescu has shown in \cite{Po89b} that for completely non-coisometric tuples,
in particular for $*$-stable ones, his characteristic function is a complete invariant
for unitary equivalence. In Section 6
we prove that our extended characteristic function does the same for the tuples $\uA$
described above. In this sense it is {\it characteristic}.
This is remarkable because the strength of Popescu's definition lies in the
completely non-coisometric situation while we always deal with a coisometric tuple $\uA$.
The extended characteristic function also does not depend on the choice of the
decomposition $\sum^d_{i=1} A_i \cdot A^*_i$ of the completely positive map $Z$
and hence also characterizes $Z$ up to conjugacy. We think that together with its nice
properties established earlier this clearly indicates that the extended
characteristic function is a valuable tool for classifying and
investigating such tuples respectively such completely positive maps.

Section 7 contains a worked example for the constructions in this paper.

\section{Weak Markov dilations and conjugacy}

In this section we give a brief and condensed review of results in
\cite{Go04}, Chapter 2, which will be used in the following and which,
as described in the introduction, motivated the investigations
documented in this paper. We also introduce notation.

A theory of {\em weak Markov dilations} has been developed in \cite{BP94}.
For a (single) normal unital completely positive map
$Z: B(\CH) \rightarrow B(\CH)$, where $B(\CH)$ consists of the bounded linear
operators on a (complex, separable) Hilbert space, it asks for a
normal unital $^*-$endomorphism
$\hat{J}: B(\hat{\CH}) \rightarrow B(\hat{\CH})$, where
$\hat{\CH}$ is a Hilbert space containing $\CH$, such that for
all $n\in\N$ and all $x\in B(\CH)$
\[
Z^n(x) = p_{\CH}\, \hat{J}^n(x p_{\CH})\; |_{\CH}.
\]
Here $p_{\CH}$ is the orthogonal projection onto $\CH$. There are
many ways to construct $\hat{J}$. In \cite{Go04}, 2.3, we gave a construction
analogous to the idea of `coupling to a shift' used in \cite{Ku85} for
describing quantum Markov processes. This gives rise to a number of
interesting problems which remain hidden in other constructions.
\\

We proceed in two steps. First note that there is a Kraus decomposition
$Z(x) = \sum^d_{i=1} a_i\, x\, a^*_i$ with $(a_i)^d_{i=1} \subset B(\CH)$.
Here $d=\infty$ is allowed in which case the sum should be interpreted
as a limit in the strong operator topology.
Let $\CP$ be a $d$-dimensional Hilbert space with
orthonormal basis $\{\epsilon_1,\ldots,\epsilon_d\}$, further $\CK$
another Hilbert space with a distinguished unit vector $\OK \in \CK$.
We identify $\CH$ with $\CH \otimes \OK \subset \CH \otimes \CK$ and
again denote by $p_{\CH}$ the orthogonal projection onto $\CH$.
For $\CK$ large enough there exists an isometry
\[
     u: \CH \otimes \CP \rightarrow \CH \otimes \CK
\quad \mbox{s.t.}\quad p_{\CH}\, u (h \otimes \epsilon_i) = a_i(h),
\]
for all $h\in\CH,\; i=1,\ldots,d$, or equivalently,
\[
u^* (h \otimes \OK) = \sum^d_{i=1} a^*_i(h) \otimes \epsilon_i.
\]
Explicitly, one may take $\CK = \C^{d+1}$ (resp. infinite-dimensional)
and identify
\[
\CH \otimes \CK \;\simeq\; (\CH \otimes \OK) \oplus \bigoplus^d_1 \CH
\;\simeq\; \CH \oplus \bigoplus^d_1 \CH.
\]

Then, using isometries
$u_1,\ldots,u_d: \CH \rightarrow \CH \oplus \bigoplus^d_1 \CH$ with
orthogonal ranges and such that $a_i = p_{\CH} u_i$ for all $i$ (for example, such
isometries are explicitly constructed in Popescu's formula
for isometric dilations, cf. \cite{Po89a} or equation 3.2 in Section 3),
we can define
\[
u (h \otimes \epsilon_i) := u_i(h)
\]
for all $h\in\CH,\; i=1,\ldots,d$ and check that $u$ has the desired
properties.
Now we define a $^*-$homomorphism
\begin{eqnarray*}
J: B(\CH) &\rightarrow& B(\CH \otimes \CK), \\
     x &\mapsto& u \, (x\otimes \eins_{\CP})\, u^*.
\end{eqnarray*}
It satisfies
\[
 p_{\CH}\, J(x) (h \otimes \OK) = p_{\CH}\,u \, (x\otimes \eins) u^*
(h \otimes \OK)
\]
\[
= p_{\CH}\,u \, (x\otimes \eins)
\big( \sum^d_{i=1} a^*_i(h) \otimes \epsilon_i \big)
= \sum^d_{i=1} a_i\, x\, a^*_i (h) = Z(x)(h),
\]
which means that $J$ is a kind of first order dilation for $Z$.
\\

For the second step we write $\tilde{\CK} := \bigotimes^\infty_1 \CK$
for an infinite tensor product of Hilbert spaces along the sequence
$(\OK)$ of unit vectors in the copies of $\CK$.
We have a distinguished unit vector $\Omega_{\tilde{\CK}}$ and
a (kind of) tensor shift
\[
R: B(\tilde{\CK}) \rightarrow B(\CP \otimes \tilde{\CK}),
\quad \tilde{y} \mapsto \eins_{\CP} \otimes \tilde{y}.
\]
Finally $\tilde{\CH} := \CH \otimes \tilde{\CK}$ and we define
a normal $^*-$endomorphism
\begin{eqnarray*}
\tilde{J}: B(\tilde{\CH}) &\rightarrow& B(\tilde{\CH}), \\
B(\CH) \otimes B(\tilde{\CK}) \ni x \otimes \tilde{y}
&\mapsto& J(x) \otimes \tilde{y} \in B(\CH \otimes \CK) \otimes B(\tilde{\CK}).
\end{eqnarray*}
Here we used von Neumann tensor products and (on the right hand side)
a shift identification $\CK \otimes \tilde{\CK} \simeq \tilde{\CK}$.
We can also write $\tilde{J}$ in the form
\[
\tilde{J}(\cdot) = u\, (Id_{\CH} \otimes R)(\cdot) \, u^*,
\]
where $u$ is identified with $u \otimes \eins_{\tilde{\CK}}$.
The natural embedding
$\CH \simeq \CH \otimes \Omega_{\tilde{\CK}} \subset \tilde{\CH}$
leads to the restriction
$\hat{J} := \tilde{J} |_{\hat{\CH}}$ with
$\hat{\CH} := \overline{\mbox{span}}_{n\geq 0} \tilde{J}^n(p_{\CH})(\tilde{\CH})$,
which can be checked to be a normal unital $^*$-endomorphism
satisfying all the properties of a weak Markov dilation for $Z$
described above. See \cite{Go04}, 2.3.

A Kraus decomposition of $\hat{J}$ can be written as
\[
\hat{J}(x) = \sum^d_{i=1} t_i\, x\, t^*_i,
\]
where $t_i \in B(\hat{\CH})$ is obtained by linear extension
of $\CH \otimes \tilde{\CK} \ni h \otimes \tilde{k}
\mapsto u_i(h) \otimes \tilde{k} = u (h \otimes \epsilon_i) \otimes \tilde{k}
\in (\CH \otimes \CK)\otimes \tilde{\CK} \simeq \CH \otimes \tilde{\CK}$.
Because $\hat{J}$ is a normal unital $^*-$endomorphism the $(t_i)^d_{i=1}$
generate a representation of the Cuntz algebra $\Od$ on $\hat{\CH}$
which we called a {\em coupling representation} in \cite{Go04}, 2.4.
Note that the tuple $(t_1,\ldots,t_d)$ is an isometric dilation of
the tuple $(a_1,\ldots,a_d)$, i.e., the $t_i$ are isometries with
orthogonal ranges and $p_{\CH} t^n_i |_{\CH} = a^n_i$ for all $i=1,\ldots,d$
and $n\in\N$.
\\

The following {\em multi-index notation} will be used frequently in this work.
Let $\Lambda $ denote the set
$\{ 1, 2, \ldots , d\}.$ For  operator tuples
$(a_1, \ldots , a_d),$ given $\alpha =(\alpha _1, \ldots , \alpha _m)$
in $\Lambda ^m$,\, $a_{\alpha }$ will stand for the operator
$a_{\alpha _1}a_{\alpha _2}\ldots a_{\alpha _m}$, \; $|\alpha| := m$.
Further $\tilde{\Lambda}:=\cup_{n=0}^{\infty} \Lambda^n$, where
$\Lambda ^0:=\{ 0\}$ and $a_0$ is the identity operator.
If we write $a^*_{\alpha}$ this always means
$(a_{\alpha})^* = a^*_{\alpha _m}\ldots a^*_{\alpha _1}$.
\\

Back to our isometric dilation, it can be checked that
\[
\overline{\mbox{span}}\{t_{\alpha} h: h \in \CH, \alpha \in \tilde{\Lambda} \}
=\hat{\CH},
\]
which means that we have a {\em minimal isometric dilation},
cf. \cite{Po89a} or the beginning of Section 3. For more details on the
construction above see \cite{Go04}, 2.3 and 2.4.

Assume now that there is an invariant vector state for
$Z: B(\CH) \rightarrow B(\CH)$ given by a unit vector $\OH \in \CH$.
Equivalent: There is a unit vector
$\OP = \sum^d_{i=1} \overline{\omega}_i \epsilon_i \in \CP$
such that $u (\OH\otimes\OP) = \OH \otimes \OK$. Also equivalent:
For $i=1,\ldots,d$ we have $a^*_i\, \OH = \overline{\omega}_i\, \OH$.
Here $\omega_i \in \C$ with $\sum^d_{i=1} |\omega_i|^2 = 1$ and
we used complex conjugation to get nice formulas later. See
\cite{Go04}, A.5.1, for a proof of the equivalences.

On $\tilde{\CP} := \bigotimes^\infty_1 \CP$ along the unit vectors
$(\OP)$ in the copies of $\CP$ we have a tensor shift
\[
S: B(\tilde{\CP}) \rightarrow B(\tilde{\CP}),
\quad \tilde{y} \mapsto \eins_{\CP} \otimes \tilde{y}.
\]
Its Kraus decomposition is $S(\tilde{y}) = \sum^d_{i=1} s_i\,\tilde{y}\,s^*_i$
with $s_i \in B(\tilde{\CP})$ and
$s_i(\tilde{k}) = \epsilon_i \otimes \tilde{k}$
for $\tilde{k}\in\tilde{\CP}$ and $i=1,\ldots,d$.
In \cite{Go04}, 2.5,
we obtained an interesting description of the situation when the
dilation $\hat{J}$ is conjugate to the shift endomorphism $S$. This
result will be further analyzed in this paper. We give a version
suitable for our present needs but the reader should have no problems to
obtain a proof of the following from \cite{Go04}, 2.5.

\begin{thm}
Let $Z: B(\CH) \rightarrow B(\CH)$ be a normal unital completely positive
map with an invariant vector state $\langle \OH, \cdot\, \OH \rangle$.
Notation as introduced above, $d\geq 2$. The following assertions are
equivalent:
\begin{itemize}
\item[(a)]
$Z$ is ergodic, i.e., the fixed point space of $Z$ consists of
multiples of the identity.
\item[(b)]
The vector state $\langle \OH, \cdot\, \OH \rangle$ is absorbing
for $Z$, i.e., if $n\to\infty$ then
$\phi(Z^n(x)) \to \langle \OH, x \OH \rangle\;$
for all normal states $\phi$ and all $x\in B(\CH)$.
(In particular, the invariant vector state is unique.)
\item[(c)]
$\hat{J}$ and $S$ are conjugate, i.e., there exists a unitary
${\bf w}: \hat{\CH} \rightarrow \tilde{\CP}$ such that
\[
\hat{J}(\hat{x}) = {\bf w}^*\, S({\bf w}\, \hat{x}\, {\bf w}^*)\, {\bf w}.
\]
\item[(d)]
The $\Od-$representations corresponding to $\hat{J}$ and $S$
are unitarily equi\-valent, i.e.,
\[
{\bf w}\, t_i = s_i\, {\bf w} \quad \mbox{for}\; i=1,\ldots,d.
\]
\end{itemize}
An explicit formula can be given for an intertwining unitary as
occurring in (c) and (d). If any of the assertions above is valid then the
 following limit exists strongly,
\[
\tilde{{\bf w}} = \lim_{n\to\infty} u^*_{0n} \ldots u^*_{01}: \;
\CH \otimes \tilde{\CK} \rightarrow \CH \otimes \tilde{\CP},
\]
where we used a leg notation, i.e.,
$u_{0n} = (Id_{\CH} \otimes R)^{n-1}(u)$.
In other words $u_{0n}$ is $u$ acting on $\CH$ and on the $n-$th copy
of $\CP$. Further $\tilde{{\bf w}}$ is a partial isometry with initial space
$\hat{\CH}$ and final space $\tilde{\CP} \simeq \OH \otimes \tilde{\CP}
\subset \CH \otimes \tilde{\CP}$ and we can define ${\bf w}$ as the
corresponding restriction of $\tilde{{\bf w}}$.
\end{thm}

To illustrate the product formula for ${\bf w}$,
which will be our main interest in this work,
we use it to derive (d).
\[
{\bf w}\,t_i (h \otimes \tilde{k}) = {\bf w} \,
\big[ u (h \otimes \epsilon_i) \otimes \tilde{k} \big]
= \lim_{n\to\infty} u^*_{0n} \ldots u^*_{01} u_{01}
(h \otimes \epsilon_i \otimes \tilde{k})
\]
\[
= \lim_{n\to\infty} u^*_{0n} \ldots u^*_{02}
(h \otimes \epsilon_i \otimes \tilde{k})
= s_i\, {\bf w} (h \otimes \tilde{k}).
\]
Let us finally note that Theorem 1.1 is related to
the conjugacy results in \cite{Pow88} and \cite{BJP96}.
Compare also Proposition 2.4.

\section{Ergodic coisometric row contractions}

In the previous section we considered a map
$Z: B(\CH) \rightarrow B(\CH)$ given by
$Z(x) = \sum^d_{i=1} A_i\, x\, A^*_i$, where $(A_i)^d_{i=1} \subset B(\CH)$.
We can think of $(A_i)^d_{i=1}$ as a $d$-tuple
$\underline{A} = (A_1,\ldots,A_d)$ or
(with the same notation) as a linear map
\[
\underline{A} = (A_1,\ldots,A_d):\; \bigoplus^d_{i=1} \CH \rightarrow \CH.
\]
(Concentrating now on the tuple we have changed to capital letters $A$.
We will sometimes return to lower case letters $a$ when we want to emphasize
that we are in the (tensor product) setting of Section 1.) We have the
following dictionary.

\begin{eqnarray*}
Z(\eins) \leq \eins
&\Leftrightarrow&
\sum^d_{i=1} A_i\,A^*_i \leq \eins \\
&\Leftrightarrow&
\underline{A} \;\mbox{is a contraction} \\
\\
Z(\eins) = \eins &
\Leftrightarrow&
\sum^d_{i=1} A_i\,A^*_i = \eins \\
\big( Z \;\mbox{is called unital}\big)
& &
\big( \underline{A} \;\mbox{is called coisometric} \big) \\
\\
\langle \OH, \cdot \OH \rangle = \langle \OH, Z(\cdot) \OH \rangle
&\Leftrightarrow&
A^*_i\, \OH = \overline{\omega}_i\, \OH, \;
\omega_i \in \C, \; \sum^d_{i=1} |\omega_i|^2 = 1 \\
\big(\mbox{ invariant vector state} \big)
& &
\big(\mbox{ common eigenvector for adjoints} \big) \\
\\
Z \;\mbox{ergodic}
&\Rightarrow&
\{A_i, A^*_i\}^\prime = \C\,\eins \\
\big(\mbox{trivial fixed point space} \big)
& &
\big(\;\mbox{trivial commutant}\; \big)
\end{eqnarray*}

The converse of the implication at the end of the dictionary
is not valid. This is related to the fact that the fixed point space of
a completely positive map is not always an algebra. Compare the detailed
discussion of this phenomenon in \cite{BJKW00}.

By a slight abuse of language we call the tuple (or row contraction)
$\underline{A} = (A_1,\ldots,A_d)$ {\em ergodic} if the corresponding map $Z$
is ergodic. With this terminology we can interpret Theorem 1.1
as a result about ergodic coisometric row contractions $\underline{A}$
with a common eigenvector $\OH$ for the adjoints $A^*_i$. This will be
examined starting with Section 3. To represent these objects more
explicitly let us write $\Ho := \CH \ominus \C\,\OH$. With respect to
the decomposition $\CH = \C\,\OH \oplus \Ho$ we get
$2 \times 2 -$ block matrices
\begin{equation}
A_i \left( \begin{array}{cc}
        \omega_i & 0 \\
        |\ell_i\rangle & \Ao_i \\
        \end{array}
\right),
\quad
A^*_i \left( \begin{array}{cc}
        \overline{\omega}_i & \langle \ell_i| \\
        0 & \Aso_i \\
        \end{array}
\right).
\end{equation}

Here $\Ao_i \in B(\Ho)$ and $\ell_i\in\Ho$. For the off-diagonal terms
we used a Dirac notation that should be clear without further comments.

Note that the case $d=1$ is rather uninteresting in this setting because
if $A$ is a coisometry with block matrix
$
\left( \begin{array}{cc}
        \omega & 0 \\
        |\ell\rangle & \Ao \\
        \end{array}
\right)
$
then because
\[
\left( \begin{array}{cc}
        1 & 0 \\
        0 & \eins \\
        \end{array}
\right)
= A\,A^* \left( \begin{array}{cc}
        |\omega|^2 & \omega\, \langle\ell| \\
        \overline{\omega}\, |\ell\rangle & |\ell\rangle\langle\ell|
          + \Ao\,\Aso\\
        \end{array}
\right)
\]
we always have $\ell = 0$. But for $d \geq 2$ there are many interesting
examples arising from unital ergodic completely positive maps with
invariant vector states. See Section 1 and also Section 7 for an
explicit example. We always assume $d \geq 2$.

\begin{prop}
A coisometric row contraction $\underline{A} = (A_1,\ldots,A_d)$ is ergodic
with common eigenvector $\OH$ for the adjoints
$A^*_1,\ldots,A^*_d$ if and only if $\Ho$ is invariant for
$A_1,\ldots,A_d$ and the restricted row contraction
$(\Ao_1,\ldots,\Ao_d)$ on $\Ho$ is $*$-stable, i.e., for all $h\in\Ho$
\[
\lim_{n\to\infty} \sum_{|\alpha|=n} \| \Aso_\alpha h\|^2 = 0\;.
\]
\end{prop}

Here we used the multi-index notation introduced in Section 1. Note that
$*$-stable tuples are also called pure, we prefer the terminology from
\cite{FF90}.
\\

\begin{proof}
It is clear that $\OH$ is a common eigenvector for the adjoints if and
only if $\Ho$ is invariant for $A_1,\ldots,A_d$. Let
$Z(\cdot) = \sum^d_{i=1} A_i \cdot A^*_i$
be the associated completely positive map. With
$q \, := \eins - | \OH \rangle \langle \OH |$, the orthogonal projection
onto $\Ho$, and by using $ q\, A_i\,q = A_i\, q \simeq\; \Ao_i$ for
all $i$, we get
\[
Z^n(q) = \sum_{|\alpha|=n} A_\alpha\, q\, A^*_\alpha
= \sum_{|\alpha|=n} \Ao_\alpha\, \Aso_\alpha
\]
and thus for all $h\in\Ho$
\[
\sum_{|\alpha|=n} \| \Aso_\alpha h\|^2 = \langle h, Z^n(q)\, h \rangle.
\]
Now it is well known that ergodicity of $Z$ is equivalent to
$Z^n(q) \to 0$ for $n\to\infty$ in the weak operator topology.
See \cite{GKL06}, Prop. 3.2.
This completes the proof.
\end{proof}

\begin{rem}
Given a coisometric row contraction $\underline{a} = (a_1,\ldots,a_d)$
we also have the isometry $u: \CH \otimes \CP \rightarrow \CH \otimes \CK$
from Section 1. We introduce the linear
map $a: \CP \rightarrow B(\CH),\; k \mapsto a_k$ defined by
\[
a^*_k(h) \otimes k := (\eins_{\CH} \otimes | k \rangle \langle k |)\,
u^*(h \otimes \OK).
\]
Compare \cite{Go04}, A.3.3. In particular $a_i = a_{\epsilon_i}$ for
$i=1,\ldots,d$, where $\{\epsilon_1,\ldots,\epsilon_d\}$
is the orthonormal basis of $\CP$ used in the definition of $u$.
Arveson's metric operator spaces, cf. \cite{Ar03}, give a conceptual
foundation for basis transformations in the operator space linearly
spanned by the $a_i$. Similarly, in our formalism a unitary in $B(\CP)$
transforms $\underline{a} = (a_1,\ldots,a_d)$ into another tuple
$\underline{a}^\prime = (a^\prime_1,\ldots,a^\prime_d)$. If $\OH$
is a common eigenvector for the adjoints $a^*_i$ then $\OH$ is also
a common eigenvector for the adjoints $(a^\prime_i)^*$ but of course
the eigenvalues are transformed to another tuple $\uomega^\prime
= (\omega^\prime_1,\ldots,\omega^\prime_d)$. We should consider the
tuples $\ua$ and $\underline{a}^\prime$ to be essentially the same.
This also means that the complex numbers $\omega_i$ are not particularly
important and they should not play a role in classification. They just
reflect a certain choice of orthonormal basis in the relevant metric
operator space. Independent of basis transformations is the vector
$\OP = \sum^d_{i=1} \overline{\omega}_i\, \epsilon_i \in \CP$
satisfying $u (\OH\otimes\OP) = \OH \otimes \OK$ (see Section 1)
and the operator $a_{\OP} = \sum^d_{i=1} \overline{\omega}_i\, a_i$.
\end{rem}
For later use we show

\begin{prop}
Let $\underline{A} = (A_1,\ldots,A_d)$ be an ergodic
coisometric row contraction such that $A^*_i\, \OH = \overline{\omega}_i\, \OH$
for all $i$, further $A_{\OP} := \sum^d_{i=1} \overline{\omega}_i\, A_i$.
Then for $n\to\infty$ in the strong operator topology
\[
(A^*_{\OP})^n \to | \OH \rangle \langle \OH |.
\]
\end{prop}

\begin{proof}
We use the setting of Section 1 to be able to apply Theorem 1.1.
From $u^* (h \otimes \OK) = \sum^d_{i=1} a^*_i(h) \otimes \epsilon_i$
we obtain
\[
u^* (h \otimes \OK) = a^*_{\OP}(h) \otimes \OP\; \oplus h^\prime\]
with $h^\prime \in \CH \otimes \OP^\perp$. Assume that $h \in \Ho$.
Because $u^*$ is isometric on $\CH \otimes \OK$ we conclude that
\begin{equation}
u^*(\OH\otimes\OK) = \OH \otimes \OP \perp u^*(h \otimes\OK)
\end{equation}
and thus also $a^*_{\OP}(h) \in \Ho$. In other words,
\[
a^*_{\OP}(\Ho) \subset\, \Ho .
\]
Let $q_n$ be the orthogonal projection from $\CH \otimes \bigotimes^n_1 \CP$
onto $\OH \otimes \bigotimes^n_1 \CP$. From Theorem 1.1 it follows that
\[
(\eins - q_n) u^*_{0n} \ldots u^*_{01} (h \otimes \bigotimes^n_1 \OK) \to 0
\quad (n\to\infty).
\]
On the other hand, by iterating the formula from the beginning,
\[
u^*_{0n} \ldots u^*_{01} (h \otimes \bigotimes^n_1 \OK)
= \big((a^*_{\OP})^n(h) \otimes \bigotimes^n_1 \OP \big) \oplus h^\prime
\]
with $h^\prime \in \CH \otimes (\bigotimes^n_1 \OP)^\perp$. It follows
that also
\[
(\eins - q_n) \big( (a^*_{\OP})^n(h) \otimes \bigotimes^n_1 \OP \big)
\to 0.
\]
But from $a^*_{\OP}(\Ho) \subset\, \Ho $ we have
$q_n \big( (a^*_{\OP})^n(h) \otimes \bigotimes^n_1 \OP \big) = 0$
for all $n$. We conclude that $(a^*_{\OP})^n(h) \to 0$ for $n\to\infty$.
Further
\[
a^*_{\OP} \OH = \sum^d_{i=1} \omega_i\, a^*_i\, \OH
= \sum^d_{i=1} \omega_i\, \overline{\omega}_i\, \OH = \OH,
\]
and the proposition is proved.
\end{proof}

The following proposition summarizes some well known properties of minimal isometric
dilations and associated Cuntz algebra representations.

\begin{prop}
Suppose $\uA$ is a coisometric tuple on $\CH$ and $\uV$ is its minimal isometric dilation.
Assume $\Omega_\CH$ is a distinguished unit vector in $\CH$ and
$\uomega = (\omega_1,\ldots,\omega_d) \in \C^d,\;\sum_i |\omega_i|^2=1$.
Then the following are equivalent.
\begin{enumerate}
\item $\uA$ is ergodic and $A^*_i\, \OH = \overline{\omega}_i\, \OH$ for all $i$.

\item $\uV$ is ergodic and
$V^*_i\, \OH = \overline{\omega}_i\, \OH$ for all $i$.

\item $V^*_i\, \OH = \overline{\omega}_i\, \OH$ and $\uV$
generates the GNS-representation of the Cuntz algebra $\CO_d = C^*\{g_1,\cdots,g_d \}$
($g_i$ its abstract generators) with respect to the Cuntz state
which maps
\[ g_\alpha\, g^*_\beta \mapsto
\omega_\alpha\, \oomega_\beta, ~~\forall \alpha,
\beta \in \tilde{\Lambda}.
\]
\end{enumerate}
Cuntz states are pure and the corresponding GNS-representations are irreducible.
\end{prop}

This Proposition clearly follows from Theorem 5.1 of \cite{BJKW00}, Theorem 3.3
and Theorem 4.1 of \cite{BJP96}.
Note that in Theorem 1.1(d) we already saw a concrete version of the corresponding
Cuntz algebra representation.

\section{A new characteristic function}

First we recall some more details of the theory of minimal isometric dilations for row contractions
(cf. \cite{Po89a}) and introduce further notation.

The full Fock space over $\C^d$ ($d\geq 2$) denoted by $\Gamma (\C^d)$ is
$$\Gamma (\C^d) := \mathbb{C}\oplus \C^d \oplus (\C^d)^{\otimes ^2}\oplus \cdots
\oplus (\C^d) ^{\otimes ^m}\oplus \cdots.
$$
$1\oplus 0\oplus \cdots$ is called the vacuum vector. Let
$\{e_1, \ldots , e_d\}$ be the standard
orthonormal basis of $\C^d$. Recall that we include $d=\infty$
in which case $\C^d$ stands for a complex separable Hilbert space
of infinite dimension. For $\alpha \in {\tilde
\Lambda}$, $e_{\alpha }$ will denote the vector $e_{\alpha
_1}\otimes e_{\alpha _2}\otimes \cdots \otimes e_{\alpha _m}$ in
the full Fock space $\Gamma (\mathbb{C}^d)$ and $e_0$ will denote
the vacuum vector. Then the (left) creation operators
$L_i$ on $\Gamma({\mathbb{C}^d})$ are defined by
\[ L_i x = e_i \otimes x \]  for $1 \leq i \leq d$ and $x
\in \Gamma({\mathbb{C}}^d).$ The row contraction
$\underline{L}= (L_1, \ldots , L_d)$ consists of isometries with
orthogonal ranges.
\\

Let $\uT=(T_1,\cdots,T_d)$ be a row contraction on a Hilbert space $\CH$.
Treating $\uT$ as a row operator from $ \bigoplus^d_{i=1} \CH$
to $\CH,$ define
$D_*:= (\eins-\uT \uT^*)^\frac{1}{2}: \CH \rightarrow \CH$ and
$D:=(\eins-\uT^*\uT)^\frac{1}{2}:
\bigoplus^d_{i=1} \CH \rightarrow \bigoplus^d_{i=1} \CH$.
This implies that
\begin{equation}
D_*= (\eins -\sum^d_{i=1} T_i T^*_i)^\frac{1}{2},
~~~D=(\delta_{ij} \eins -T_i^*T_j)^\frac{1}{2}_{d\times d}.
\end{equation}
Observe that $\uT D^2 = D^2_* \uT$ and hence $\uT D = D_* \uT.$
Let $\CD:=\mbox{Range~} D$ and $\CD_*:=\mbox{Range~} D_*.$
Popescu in \cite{Po89a} gave the following explicit presentation of the
minimal isometric dilation of $\uT$ by
$\uV$ on $\CH  \oplus (\Gamma (\C^d) \otimes \CD)$,
\begin{equation}
V_i(h \oplus \sum_{\alpha \in \tilde{\Lambda}}  e_\alpha \otimes d_\alpha)
= T_i h \oplus [e_0 \otimes D_i h + e_i \otimes \sum_{\alpha \in
\tilde{\Lambda}} e_\alpha \otimes d_\alpha]
\end{equation}
for $h \in \CH$ and $d_\alpha \in \CD.$ Here
$D_i h := D (0,\ldots,0,h,0,\ldots,0)$ and $h$ is embedded at the
$i^{th}$ component.

In other words, the $V_i$ are isometries with orthogonal ranges such that
$T^*_i = V^*_i |_{\CH}$ for $i=1,\ldots,d$ and the spaces
$V_\alpha \CH$ with $\alpha \in \tilde{\Lambda}$ together span the Hilbert space
on which the $V_i$ are defined. It is an important fact, which we shall
use repeatedly, that such minimal isometric dilations are unique up to unitary
equivalence (cf. \cite{Po89a}).
\\

Now, as in Section 2,
let $\uA =(A_1,\cdots,A_d),\, A_i \in B(\CH)$, be an ergodic coisometric tuple
with $A^*_i \Omega_\CH = \oomega_i \Omega_\CH$
for some unit vector $\Omega_\CH \in \CH$ and some
$\uomega \in \C^d,\;\sum_i |\omega_i|^2=1$.
Let $\uV = (V_1,\cdots,V_d)$ be the minimal isometric dilation of $\uA$
given by Popescu's construction (see equation 3.2) on
$\CH \oplus \big(\Gamma (\C^d) \otimes \CD_A \big)$.
Because $A^*_i = V^*_i |_\CH$ we also have $V^*_i \Omega_\CH = \oomega_i \Omega_\CH$
and because $\uV$ generates an irreducible $\Od-$representation (Proposition 2.4),
we see that $\uV$ is also a minimal isometric dilation of
$\uomega: \C^d \rightarrow \C$. In fact, we can think of $\uomega$ as the
most elementary example of a tuple with all the properties stated for $\uA$.
Let $\tilde{\uV} = (\tilde{V}_1,\cdots,\tilde{V}_d)$ be the minimal isometric
dilation of $\uomega$ given by Popescu's construction on
$\C \oplus (\Gamma (\C^d) \otimes \CD_\omega)$.
\\

Because $\uA$ is coisometric it follows from equation 3.1 that $D$ is in fact
a projection and hence $D=(\delta_{ij}\eins-A_i^*A_j)_{d\times d}.$
We infer that $D (A^*_1, \cdots,A^*_d)^T =0$, where $T$ stands for
transpose.
Applied to $\uomega$
instead of $\uA$ this shows that $D_{\omega} =
(\eins - |\overline{\uomega} \rangle \langle \overline{\uomega} |)$ and
$$\CD_\omega \oplus \C (\oomega_1, \cdots, \oomega_d)^T = \C^d,$$
where $\underline{\oomega} = (\oomega_1, \cdots, \oomega_d)$.
\\

\begin{rem}
Because $\OH$ is cyclic for $\{ V_\alpha,\, \alpha \in \tilde{\Lambda} \}$
we have
$$\overline{\mbox{span}}\{A_\alpha \Omega_{\CH}:\alpha \in
\tilde{\Lambda}\}=\overline{\mbox{span}}\{p_\CH\,V_\alpha \Omega_{\CH}:
\alpha \in \tilde{\Lambda}
\}=\; \CH.$$
Using the notation from equation 2.1 this further implies that
$$\overline{\mbox{span}}\{ \Ao_\alpha\, l_i: \alpha \in \tilde{\Lambda},
1 \leq i \leq d\} =\; \Ho .$$
\end{rem}
\vspace{0.4cm}

As minimal isometric dilations of the tuple $\uomega$ are unique up
to unitary equi\-valence, there exists a unitary
$$W: \CH \oplus (\Gamma (\C^d) \otimes \CD_A) \to
\C \oplus (\Gamma (\C^d) \otimes \CD_\omega),$$
such that $W V_i = \tilde{V}_i W$ for all $i.$

After showing the existence of $W$ we now proceed to compute $W$ explicitly.
For $\uA$, by using Popescu's construction, we have its minimal
isometric dilation $\uV$ on $\CH \oplus (\Gamma (\C^d) \otimes \CD_A).$
Another way of constructing a minimal isometric dilation $\ut$ of $\ua$
was demonstrated in
Section 1 on the space $\hat{\CH}$ (obtained by restricting to the minimal subspace
of $\CH \otimes \tilde{\CK}$ with respect to $\ut$).
Identifying $\uA$ and $\ua$ on the Hilbert space $\CH$ there is a unitary
$\Gamma_A: \hat{\CH} \rightarrow \CH \oplus (\Gamma (\C^d) \otimes \CD_A)$
which is the identity on $\CH$ and satisfies $V_i \Gamma_A = \Gamma_A t_i$.
\\

By Theorem 1.1(d) the tuple $\us$ on $\tilde{\CP}$ arising from the
tensor shift is unitarily equivalent to $\ut$ (resp. $\uV$), explicitly
${\bf w}\,t_i = s_i\,{\bf w}$ for all $i$. An alternative viewpoint
on the existence of ${\bf w}$
is to note that $\us$ is a minimal isometric dilation of $\uomega.$ In fact,
$s^*_i\, \Omega_{\tilde{\CP}} = \langle \epsilon_i, \OP \rangle \Omega_{\tilde{\CP}}
= \oomega_i\, \Omega_{\tilde{\CP}}$ for all $i$. Hence there is also a unitary
$\Gamma_\omega: \tilde{\CP} \rightarrow \C \oplus (\Gamma (\C^d)
\otimes \CD_\omega)$
with $\Gamma_\omega \Omega_{\tilde{\CP}} = 1 \in \C$ which satisfies
$\tilde{V}_i \Gamma_\omega = \Gamma_\omega s_i$.
\\

\begin{rem}
It is possible to describe $\Gamma_\omega$ in an explicit way and
in doing so to construct an interesting and natural (unitary) identification
of $\bigotimes^\infty_1 \C^d$ and $\C \oplus (\Gamma (\C^d) \otimes \C^{d-1})$.
In fact, recall (from Section 1) that $\tilde{\CP} = \bigotimes^\infty_1 \CP$ and
the space $\CP$ is nothing but a $d$-dimensional Hilbert space. Hence we can identify
\[
\C^d \simeq \CP =\, \Po \oplus\, \C \OP \simeq
\CD_\omega \oplus \C\, \underline{\oomega}^T \simeq \C^{d-1} \oplus \C
\]
In this identification the orthonormal basis $(\epsilon_i)^d_{i=1}$ of $\CP$
goes to the canonical basis $(e_i)^d_{i=1}$ of $\C^d$, in particular the vector
$\OP = \sum_i \oomega_i\, \epsilon_i$
goes to $\underline{\oomega}^T = (\oomega_1, \cdots, \oomega_d)^T$
and we have $\Po\, \simeq \CD_\omega$.
Then we can write
\begin{eqnarray*}
\Gamma_\omega:\qquad \Omega_{\tilde{\CP}} &\mapsto& 1 \in \C, \\
k \otimes \Omega_{\tilde{\CP}} &\mapsto& e_0 \otimes k \\
\epsilon_\alpha \otimes k \otimes \Omega_{\tilde{\CP}} &\mapsto& e_\alpha \otimes k,
\end{eqnarray*}
where $k \in \Po,\, \alpha \in \tilde{\Lambda},\,
\epsilon_\alpha = \epsilon_{\alpha_1} \otimes \ldots \epsilon_{\alpha_n} \in
\bigotimes^n_1 \CP$ (the first $n$ copies of $\CP$ in the infinite tensor
product $\tilde{\CP}$),
$e_\alpha = e_{\alpha_1} \otimes \ldots e_{\alpha_n} \in \Gamma(\C^d)$ as usual.
It is easily checked that $\Gamma_\omega$ given in this way indeed satisfies
the equation $\tilde{V}_i \Gamma_\omega = \Gamma_\omega s_i$ (for all $i$), which may thus be seen
as the abstract characterization of this unitary map (together with $\Gamma_\omega \Omega_{\tilde{\CP}} = 1$).
\end{rem}
\vspace{0.4cm}

Summarizing, for $i=1,\ldots,d$
\begin{eqnarray*}
V_i\, \Gamma_A = \Gamma_A\, t_i,\quad
{\bf w}\, t_i = s_i\, {\bf w},\quad
\tilde{V}_i\, \Gamma_\omega = \Gamma_\omega\, s_i
\end{eqnarray*}
and we have the commuting diagram
\begin{equation}
\xymatrix{
\hat{\CH} \ar[r]^{{\bf w}} \ar[d]_{\Gamma_A}
& \tilde{\CP} \ar[d]^{\Gamma_\omega}
\\
\CH \oplus (\Gamma (\C^d) \otimes \CD_A) \ar[r]^{W}
& \C \oplus (\Gamma (\C^d) \otimes \CD_\omega).
}
\end{equation}
From the diagram we get
$$W=\Gamma_\omega {\bf w} \Gamma^{-1}_A.$$
Combined with the equations above this yields $W V_i = \tilde{V}_i\,W$ and
we see that $W$ is nothing but the dilations-intertwining map which
we have already introduced earlier. Hence ${\bf w}$ and $W$ are essentially the
same thing and for the study of certain problems it may be helpful to switch
from one picture to the other.

In the following we analyze $W$ to arrive at an interpretation as a new kind of
characteristic function. First we have an isometric embedding
\begin{equation}
\hat{C} := W |_\CH : \CH \rightarrow \C \oplus (\Gamma (\C^d) \otimes \CD_\omega).
\end{equation}
Note that $\hat{C}\,\OH = W\,\OH = 1 \in \C$. The remaining part is an isometry
\begin{equation}
M_{\hat{\Theta}} := W |_{\Gamma (\C^d) \otimes \CD_A} :
\Gamma (\C^d) \otimes \CD_A \rightarrow \Gamma (\C^d) \otimes \CD_\omega .
\end{equation}
From equation 3.2 we get for all $i$
\[
V_i\, |_{\Gamma (\C^d) \otimes \CD_A} = (L_i \otimes \eins_{\CD_A}),
\]
\[
\tilde{V}_i\, |_{\Gamma (\C^d) \otimes \CD_\omega} = (L_i \otimes \eins_{\CD_\omega}),
\]
and we conclude that
\begin{equation}
M_{\hat{\Theta}} (L_i \otimes \eins_{\CD_A}) =
(L_i \otimes \eins_{\CD_\omega}) M_{\hat{\Theta}},~~\forall 1 \leq i \leq d.
\end{equation}
In other words, $M_{\hat{\Theta}}$ is a multi-analytic inner function in the sense of
\cite{Po89c,Po95}. It is determined by its symbol
\begin{equation}
\hat{\theta} := W|_{e_0 \otimes \CD_A} : \CD_A \rightarrow \Gamma (\C^d) \otimes \CD_\omega,
\end{equation}
where we have identified $e_0 \otimes \CD_A$ and $\CD_A$. In other words, we think of the
symbol $\hat{\theta}$ as an isometric embedding of $\CD_A$ into $\Gamma (\C^d) \otimes \CD_\omega$.

\begin{defn}
We call $M_{\hat{\Theta}}$ (or $\hat{\theta}$) the extended characteristic function of the
row contraction $\uA$,
\end{defn}

See Sections 5 and 6 for more explanation and justification of this terminology.

\section{Explicit computation of the extended characteristic function}

To express the extended characteristic function more explicitly in terms of the tuple $\uA$
we start by defining
\begin{eqnarray}
\hat{D}_*:\; \Ho = \CH \ominus \C \OH &\rightarrow& \Po\,
= \CP \ominus \C \OP \simeq \CD_\omega,
\end{eqnarray}
\[
h \mapsto \big( \langle \OH | \otimes \eins_{\CP} \big) \, u^*(h \otimes \OK),
\]
where $u: \CH \otimes \CP \rightarrow \CH \otimes \CK$ is the isometry introduced
in Section 1. That indeed the range of $\hat{D}_*$ is contained in $\Po$ follows
from equation 2.2, i.e., $u^*(h \otimes\OK) \perp \OH \otimes \OP$ for $h \in \Ho$.
With notations from equation 2.1 we can get a more concrete formula.

\begin{lem}
For all $h \in \Ho$ we have
$\hat{D}_*(h) = \sum^d_{i=1} \langle \ell_i, h \rangle \epsilon_i$.
\end{lem}

\begin{proof}
$
\big( \langle \OH | \otimes \eins_{\CP} \big) \, u^*(h \otimes \OK)
= \sum^d_{i=1} \langle \OH, a^*_i h \rangle \otimes \epsilon_i
= \sum^d_{i=1} \langle \ell_i, h \rangle \epsilon_i.
$
\end{proof}
\vspace{0.2cm}

\begin{prop}
The map $\hat{C}:\CH \to \C \oplus (\Gamma (\C^d) \otimes \CD_\omega)$
from equation 3.4 is given explicitly by
$\hat{C}\Omega_\CH = 1$ and for $h \in \Ho$ by
$$\hat{C} h = \sum_{\alpha\in\tilde{\Lambda}} e_\alpha \otimes \hat{D}_* \Aso_\alpha h.$$
\end{prop}
\begin{proof}
As $W \Omega_\CH = 1$ also $\hat{C} \Omega_\CH =1$.
Assume $h \in \Ho$. Then
\begin{eqnarray*}
u_{01}(h \otimes \Omega_{\tilde{\CK}}) &=&\sum_i a^*_i h \otimes \epsilon_i
\otimes \Omega_{\tilde{\CK}}\\
& = & \sum_i \langle \ell_i, h \rangle \Omega_\CH \otimes \epsilon_i
\otimes \Omega_{\tilde{\CK}} + \sum_i \aso_i\!h \otimes \epsilon_i \otimes
\Omega_{\tilde{\CK}}.
\end{eqnarray*}
Because $u^*(\OH \otimes \OK) = \OH \otimes \OP$ we obtain (with Lemma 4.1)
for the first part
\begin{eqnarray*}
&&\lim_{n \to \infty} u^*_{0n} \cdots u^*_{02}
(\sum_i \langle \ell_i, h \rangle \Omega_\CH \otimes \epsilon_i
\otimes \Omega_{\tilde{\CK}})\\
&=&\sum_i \langle \ell_i, h \rangle \Omega_\CH \otimes \epsilon_i
\otimes \Omega_{\tilde{\CP}} = \OH \otimes \hat{D}_* h \otimes \Omega_{\tilde{\CP}}
\simeq \hat{D}_* h \otimes \Omega_{\tilde{\CP}} \in \tilde{\CP}.
\end{eqnarray*}
Using the product formula from Theorem 1.1 and iterating the argument above we get
\[
\hat{C}(h) \,=\, Wh \,=\, \Gamma_\omega {\bf w} \Gamma^{-1}_A(h)
\]
\[
= \Gamma_\omega(\hat{D}_* h \otimes \Omega_{\tilde{\CP}})
\,+\, \Gamma_\omega  \lim_{n \to \infty} u^*_{0n} \cdots u^*_{02}
\sum_i \aso_i\!h \otimes \epsilon_i \otimes \Omega_{\tilde{\CK}}
\]
\[
= e_0 \otimes \hat{D}_* h
\,+\, \Gamma_\omega  \lim_{n \to \infty} u^*_{0n} \cdots u^*_{03}
\sum_{j,i} \big(\langle \ell_j, \aso_i\!h \rangle \OH \,+\, \aso_j \aso_i\!h
\big)
\otimes \epsilon_i \otimes \epsilon_j \otimes \Omega_{\tilde{\CK}}
\]
\[
= e_0 \otimes \hat{D}_* h + \sum^d_{i=1} e_i \otimes \hat{D}_* \aso_i\!h
\,+\, \Gamma_\omega  \lim_{n \to \infty} u^*_{0n} \cdots u^*_{03}
\sum_{j,i} \aso_j \aso_i\!h \otimes \epsilon_i \otimes \epsilon_j \otimes \Omega_{\tilde{K}}
\]
\[
= \ldots
\]
\[
= \sum_{|\alpha|<m} e_\alpha \otimes \hat{D}_* \aso_\alpha h
\,+\, \Gamma_\omega  \lim_{n \to \infty} u^*_{0n} \cdots u^*_{0,m+1}
\sum_{|\alpha|=m} \aso_\alpha h \otimes \epsilon_\alpha \otimes
\Omega_{\tilde{\CK}}.
\]
From Proposition 2.1 we have $\sum_{|\alpha|=m} \| \aso_\alpha\!h \|^2 \to 0$
for $m \to \infty$ and we conclude that the last term converges to $0$. It follows
that the series converges and this proves Proposition 4.2.
\end{proof}
\begin{rem}
Another way to prove Proposition 4.2 for $h \in \Ho$ consists in
repeatedly applying the formula
\[
u^* (h\otimes \OK) = a^*_{\OP} h \otimes \OP + h^\prime, \quad
h^\prime \in \CH \otimes \Po
\]
to the $u^*_{0n}(h\otimes \OK)$ and then using $(a^*_{\OP})^n h \to 0$,
see Proposition 2.3. This gives some insight how the infinite product in Theorem 1.1
transforms into the infinite sum in Proposition 4.2.
\end{rem}

Now we present an explicit computation of the extended characteristic function. One way of writing $\CD_A$ is
$$\CD_A= \overline{\mbox{span}}\{ (V_i - A_i) h: i \in \Lambda, h \in \CH\}.$$
Let $d^i_h:=(V_i -A_i)h.$ Then
$$\hat{\theta}\, d^i_h =W(V_i-A_i) h = \tilde{V}_i \hat{C} h - \hat{C} A_i h.$$

{\bf Case I:} Take $h = \Omega_\CH.$
$$\tilde{V}_i \hat{C} \Omega_\CH = \tilde{V}_i 1 = \omega_i \oplus
[e_0 \otimes (\eins -  |\underline{\oomega} \rangle \langle
\underline{\oomega} |) \epsilon_i],$$
$$\hat{C}A_i\, \Omega_\CH = \omega_i \oplus \sum_\alpha e_\alpha \otimes \hat{D}_* \Aso_\alpha l_i$$
and thus
\begin{eqnarray*}
\hat{\theta}\; d^i_{\Omega_\CH} &=& e_0 \otimes
[(\eins -  |\underline{\oomega} \rangle \langle \underline{\oomega} |)
\epsilon_i- \hat{D}_*l_i] - \sum_{|\alpha|\geq 1} e_\alpha \otimes \hat{D}_*
\Aso_\alpha l_i\\
&=& e_0 \otimes [\epsilon_i - \sum_j \oomega_j \omega_i \epsilon_j - \sum_j
\langle l_j, l_i\rangle \epsilon_j] - \sum_{|\alpha|\geq 1} e_\alpha \otimes
\sum_j \langle l_j,
\Aso_\alpha l_i \rangle \epsilon_j\\
&=&e_0 \otimes [\epsilon_i - \sum_j (\oomega_j \omega_i + \langle l_j,
l_i\rangle) \epsilon_j] - \sum_{|\alpha|\geq 1} e_\alpha \otimes \sum_j
\langle \Ao_\alpha l_j,
l_i \rangle \epsilon_j
\end{eqnarray*}
\begin{equation}
= e_0 \otimes [\epsilon_i - \sum_j
\langle A_j\,\OH, A_i\,\OH \rangle\,
\epsilon_j] - \sum_{|\alpha|\geq 1} e_\alpha \otimes \sum_j
\langle \Ao_\alpha l_j,
l_i \rangle \epsilon_j.
\end{equation}

{\bf Case II:} Now let $h \in\, \Ho$. With $i \in \Lambda$
\[
\tilde{V}_i\, \hat{C} h = (L_i \otimes \eins) \hat{C} h
= \sum_\alpha e_i \otimes e_\alpha \otimes \hat{D}_* \Aso_\alpha h,$$
$$\hat{C}A_i h= \sum_\beta e_\beta \otimes \hat{D}_* \Aso_\beta \Ao_i h.
\]
Finally
\begin{equation*}
\hat{\theta}\; d^i_h =  \sum_\alpha e_i \otimes e_\alpha \otimes \hat{D}_* \Aso_\alpha h- \sum_\beta e_\beta \otimes \hat{D}_* \Aso_\beta \Ao_i h
\end{equation*}
\begin{equation*}
=\; - e_0 \otimes \hat{D}_* \Ao _i h +  e_i \otimes \sum_\alpha e_\alpha \otimes \hat{D}_* \Aso_\alpha (\eins - \Aso_i \Ao_i) h
+ \sum_{j \neq i} e_j \otimes
\sum_\alpha e_\alpha \otimes \hat{D}_* \Aso_\alpha (-\Aso_j \Ao_i) h
\end{equation*}
\begin{equation}
=\; - e_0 \otimes \hat{D}_* \Ao _i h + \sum^d_{j=1} e_j \otimes
\sum_{\alpha \in \tilde{\Lambda}} e_\alpha \otimes \hat{D}_*
\Aso_\alpha (\delta_{ji} \eins - \Aso_j \Ao_i) h .
\end{equation}

\section{Case II is the characteristic function of $\uAo$}

In this section we show that case II in the previous section can be identified with the
characteristic function of the $*$-stable tuple $\uAo$, in the sense introduced by Popescu
in \cite{Po89b}. This is the reason why we have called $\hat{\theta}$ an {\it extended}
characteristic function. All information about $\uA$ beyond $\uAo$ must be contained in case I.

First recall the theory of characteristic functions for row contractions, as developed by G.\,Popescu
in \cite{Po89b}, generalizing the theory of B.\,Sz.-Nagy and C.\,Foias (cf. \cite{NF70}) for
single contractions. We only need the results about a $*$-stable tuple $\uAo = (\Ao_1,\ldots,\Ao_d)$
on $\Ho$. In this case, with $\Dso \;=\; (\eins- \uAo \uAso)^{\frac{1}{2}}:\; \Ho \rightarrow \Ho$ and $\CDso$ its range, the map
\begin{equation}
\Co:\; \Ho \to \Gamma (\C^d) \otimes \CDso
\end{equation}
\[
h \mapsto \sum_{\alpha \in \tilde{\Lambda}} e_\alpha \otimes \Dso \Ao_\alpha^* h
\]
is an isometry (Popescu's Poisson kernel). If, as usual,
$\Do \;=\; (\eins - \uAso \uAo)^{\frac{1}{2}}: \bigoplus^d_1 \Ho \rightarrow \bigoplus^d_1 \Ho$,
with $\CDo$ its range, and if $P_j$ is the projection onto the $j$-th component,
then the characteristic function $\theta_{\Ao}$ of $\uAo$ can be defined as
\begin{equation}
\theta_{\Ao}:\; \CDo \to \Gamma (\C^d) \otimes \CDso
\end{equation}
\[
f \mapsto -e_0 \otimes \sum^d_{j=1} \Ao_j P_j f + \sum^d_{j=1} e_j \otimes
\sum_{\alpha \in \tilde{\Lambda}} e_\alpha  \otimes \Dso \Ao_\alpha^* P_j \Do f.
\]
See \cite{Po89b} for details, in particular for the important result that
$\theta_{\Ao}$ characterizes the $*$-stable tuple $\uAo$ up to unitary equivalence.
\\

Now consider again the tuple $\uA$ of the previous section, with extended characteristic function
$\hat{\theta}$. From equation 2.1
\[
A_i \left( \begin{array}{cc}
        \omega_i & 0 \\
        |\ell_i\rangle & \Ao_i \\
        \end{array}
\right),
\quad
A^*_i \left( \begin{array}{cc}
        \overline{\omega}_i & \langle \ell_i| \\
        0 & \Aso_i \\
        \end{array}
\right)
\]
and hence
\[ A_i A^*_i = \left( \begin{array}{cc} |\oomega_i|^2 & \langle \oomega_i l_i |\\
|\oomega_i l_i \rangle & |l_i \rangle \langle l_i| + \Ao_i \Aso_i
\end{array} \right).
\]
Recall that $D^2_*= \eins - \sum_i A_i A^*_i$ which is $0$ as $\uA$
is coisometric.
Thus $\sum_i \oomega_i l_i =0$ and
$\eins-\sum_i \Ao_i \Aso_i= \sum_i |l_i \rangle \langle l_i|.$
The first equation means that $A^*_{\OP}(\Ho) \subset \Ho$
and that
\[
\langle \hat{\CD}_* h, \OP \rangle =
\langle \sum_i \langle \ell_i, h \rangle \epsilon_i,
\sum_j \oomega_j \epsilon_j \rangle
= \langle \sum_i \oomega_i \ell_i, h \rangle = 0,
\]
which we already know (see 4.1). \\
The second equation yields
\[
\Dso^2 = \eins-\sum_i \Ao_i \Aso_i= \sum_i |l_i \rangle \langle l_i|.
\]

\begin{lem}
There exists an isometry $\gamma:\, \CDso\, \rightarrow\, \Po\, \simeq \CD_{\omega}$
defined for $h \in \Ho$ as
\[
\Dso h \mapsto \sum_i \langle l_i,h \rangle \epsilon_i = \hat{D}_* h.
\]
\end{lem}
\begin{proof}
Take $h \in \Ho.$ By Lemma 4.1 we have
$\hat{D}_*(h) = \sum^d_{i=1} \langle \ell_i, h \rangle \epsilon_i$.
Now we can compute
\[
\|\hat{D}_* h\|^2 =
\langle \sum_i \langle l_i,h \rangle \epsilon_i, \sum_j \langle l_j,h \rangle
\epsilon_j\rangle = \sum_i \langle h,l_i \rangle \langle l_i,h \rangle
= \langle h, \Dso^2 h\rangle=\|\Dso h\|^2.
\]
Hence $\gamma:\; \Dso h \mapsto \hat{D}_* h$ is isometric.
\end{proof}

\begin{thm}
Let $\uA =(A_1,\cdots,A_d),\, A_i \in B(\CH)$, be an ergodic coisometric tuple
with $A^*_i \Omega_\CH = \oomega_i \Omega_\CH$
for some unit vector $\Omega_\CH \in \CH$ and some
$\uomega \in \C^d,\;\sum_i |\omega_i|^2=1$.
Let $\hat{\theta}$ be the extended characteristic function of $\uA$
and let $\theta_{\Ao}$ be the characteristic function of the ($*$-stable)
tuple $\underline{\Ao}$. For $h \in \Ho$
\begin{eqnarray*}
\gamma \Dso h &=& \hat{D}_*h,\\
(\eins \otimes \gamma)\Co h &=& \hat{C}h,
\end{eqnarray*}
\[
(\eins \otimes \gamma)\, \theta_{\Ao}\, d^i_h = \hat{\theta}\, d^i_h~~
\mbox{for}~ i \in \Lambda.
\]
In other words, the part of $\hat{\theta}$ described by case II in the previous section
is equivalent to $\theta_{\Ao}$.
\end{thm}

\begin{proof}
We only have to use Lemma 5.1 and compare Proposition 4.2 and equation 5.1 as well as
equations 4.3 and 5.2.
For the latter note that $d^i_h = \Do(0,\ldots,0,h,0\ldots,0)$, where $h$ is embedded
at the $i$-th position. Hence
\[
\gamma \sum_j \Ao_j P_j d^i_h = \gamma \sum_j \Ao_j P_j \Do(0,\ldots,0,h,0\ldots,0)
= \gamma \uAo \; \Do(0,\ldots,0,h,0\ldots,0)
\]
\[
= \gamma \Dso \uAo(0,\ldots,0,h,0\ldots,0)
= \hat{D}_* \, \Ao_i h
\]
and also
\[
P_j \Do d^i_h = P_j \Do^2(0,\ldots,0,h,0\ldots,0)
= (\delta_{ji} \eins - \Aso_j \Ao_i) h .
\]
\end{proof}

Of course, Theorem 5.2 explains why we have called $\hat{\theta}$
an {\it extended} characteristic function.

\section{The extended characteristic function is a complete unitary invariant}

In this section we prove that the extended characteristic function is a
complete invariant with respect to unitary equivalence for the row contractions
investigated in this paper. Suppose that $\uA = (A_1,\ldots,A_d)$ and
$\uB = (B_1,\ldots,B_d)$ are ergodic and coisometric row contractions on Hilbert
spaces $\CH_A$ and $\CH_B$ such that $A^*_i \Omega_A = \oomega_i \Omega_A$
and $B^*_i \Omega_B = \oomega_i \Omega_B$ for $i=1,\ldots,d$, where
$\Omega_A \in \CH_A$ and $\Omega_B \in \CH_B$ are unit vectors and
$\uomega = (\omega_1,\ldots,\omega_d)$ is a tuple of complex numbers.
Recall from Remark 2.2 that it is no serious restriction of generality
to assume that it is the same tuple of complex numbers in both cases because this can
always be achieved by a transformation with a unitary $d\times d-$matrix (with scalar entries).
We will use all the notations introduced earlier with subscripts $A$ or
$B$.

Let us say that the extended characteristic functions $\hat{\theta}_A$
and $\hat{\theta}_B$ are {\em equivalent} if there exists a unitary
$V: \CD_A \rightarrow \CD_B$ such that $\hat{\theta}_A = \hat{\theta}_B\,V$.
Note that the ranges of $\hat{\theta}_A$ and $\hat{\theta}_B$ are both
contained in $\Gamma (\C^d) \otimes \CD_\omega$ and thus this definition
makes sense. Let us further say that $\uA$ and $\uB$ are {\em unitarily
equivalent} if there exists a unitary $U: \CH_A \rightarrow \CH_B$ such
that $ U\,A_i = B_i\,U$ for $i=1,\ldots,d$. By ergodicity the unit eigenvector
$\Omega_A$ (resp. $\Omega_B$) is determined up to an unimodular constant
(see Theorem 1.1(b)) and hence in the case of unitary equivalence we can
always modify $U$ to satisfy additionally $U\,\Omega_A = \Omega_B$.

\begin{thm}
The extended characteristic functions $\hat{\theta}_A$ and $\hat{\theta}_B$
are equivalent if and only if $\uA$ and $\uB$ are unitarily equivalent.
\end{thm}
\begin{proof}
If $\uA$ and $\uB$ are unitarily equivalent then all constructions differ
only by naming and it follows that $\hat{\theta}_A$ and $\hat{\theta}_B$
are equivalent. Conversely, assume that there is a unitary
$V: \CD_A \rightarrow \CD_B$ such that $\hat{\theta}_A = \hat{\theta}_B\,V$.
Now from the commuting diagram 3.3 and the definitions following it
\begin{eqnarray*}
W_B \CH_B &=& \C \oplus \big( \Gamma (\C^d) \otimes \CD_\omega \big)
\ominus M_{\hat{\Theta}_B} \big( \Gamma (\C^d) \otimes \CD_B \big) \\
&=& \C \oplus \big( \Gamma (\C^d) \otimes \CD_\omega \big)
\ominus M_{\hat{\Theta}_B} \big( \Gamma (\C^d) \otimes V\,\CD_A \big) \\
&=& \C \oplus \big( \Gamma (\C^d) \otimes \CD_\omega \big)
\ominus M_{\hat{\Theta}_A} \big( \Gamma (\C^d) \otimes \CD_A \big) \\
&=& W_A \CH_A,
\end{eqnarray*}
where we used equation 3.6, i.e.,
$M_{\hat{\Theta}} (L_i \otimes \eins_{\CD}) =
(L_i \otimes \eins_{\CD_\omega}) M_{\hat{\Theta}},~~\forall 1 \leq i \leq d$,
to deduce $M_{\hat{\Theta}_A} = M_{\hat{\Theta}_B} (\eins \otimes V)$ from
$\hat{\theta}_A = \hat{\theta}_B\,V$.
Now we define the unitary $U$ by
\[
U := W^{-1}_B\,W_A |_{\CH_A}: \CH_A \rightarrow \CH_B.
\]
Because $W_A \Omega_A = 1 = W_B \Omega_B$ we have $U\,\Omega_A = \Omega_B$.
Further for all  $i=1,\ldots,d$ and $h\in \CH_A$,
\[
U A_i\, h = W^{-1}_B\,W_A\, A_i\,h = W^{-1}_B\,W_A P_{\CH_A}\,V^A_i h
= P_{\CH_B} W^{-1}_B\,W_A\,V^A_i h
\]
\[
= P_{\CH_B} W^{-1}_B\,\tilde{V}_i\,
W_A h = P_{\CH_B} V^B_i\,W^{-1}_B\,W_A h = B_i\,U h,
\]
i.e., $\uA$ and $\uB$ are unitarily equivalent.
\end{proof}
\begin{rem}
An analogous result for completely non-coisometric tuples has been shown by
G.\,Popescu in \cite{Po89b}, Theorem 5.4.
\end{rem}
Note further that if we change  $\uA = (A_1,\ldots,A_d)$ into
$\uA^\prime = (A^\prime_1,\ldots,A^\prime_d)$ by applying a
unitary $d\times d-$matrix with scalar entries (as described in Remark 2.2), then
$\hat{\theta}_A = \hat{\theta}_{A^\prime}$. In fact, this
follows immediately from the definition of $W$ as an intertwiner in Section 3,
from which it is evident that $W$ does not change if we take the same linear combinations on
the left and on the right. This does not contradict Theorem 6.1
because $\uomega$ and $\uomega\prime$ are now different tuples of
eigenvalues and Theorem 6.1 is only applicable when the same tuple
of eigenvalues is used for $\uA$ and $\uB$.

For another interpretation, let $Z$ be a normal, unital, ergodic, completely positive map
with an invariant vector state $\langle \Omega_A,\cdot\, \Omega_A \rangle$.
If we consider two minimal Kraus decompositions of $Z$, i.e.,
\[
Z = \sum^d_{i=1} A_i \cdot A^*_i = \sum^d_{i=1} A^\prime_i \cdot (A^\prime_i)^*,
\]
with $d$ minimal, then the tuples $\uA = (A_1,\ldots,A_d)$ into
$\uA^\prime = (A^\prime_1,\ldots,A^\prime_d)$ are related in the way considered
above (see for example \cite{Go04}, A.2). It follows that
$\hat{\theta}_A = \hat{\theta}_{A^\prime}$ does not depend on the
decomposition but can be associated to $Z$ itself. Hence we have
the following reformulation of Theorem 6.1.

\begin{cor}
Let $Z_1,\,Z_2$ be normal, unital, ergodic, completely positive maps
on $B(\CH_1),\,B(\CH_2)$ with invariant vector states
$\langle \Omega_1,\cdot\, \Omega_1 \rangle$ and
$\langle \Omega_2,\cdot\, \Omega_2 \rangle$.
Then the associated extended characteristic functions
$\hat{\theta}_1$ and $\hat{\theta}_2$ are equivalent if and only if
$Z_1$ and $Z_2$ are conjugate, i.e., there exists a unitary
$U: \CH_1 \rightarrow \CH_2$ such that
\[
Z_1(x) = U^* Z_2( U x U^* ) U \quad \mbox{for all}\; x \in B(\CH_1).
\]
\end{cor}

\section{Example}

The following example illustrates some of the constructions in this paper.
\\

Consider $\CH=\C^3$ and
$$A_1=\frac{1}{\sqrt{2}}
\left( \begin{array}{ccc} 0 & 0 & 0 \\ 1 & 0 & 0 \\ 0 & 1 & 1
\end{array} \right), A_2= \frac{1}{\sqrt{2}}
\left( \begin{array}{ccc} 1 & 1 & 0 \\ 0 &
0 & 1 \\ 0 & 0 & 0 \end{array} \right).$$
Then $\sum^2_{i=1} A_iA^*_i=\eins.$ Take the unital completely
positive map
$Z:M_3 \to M_3$ by $Z(x)=\sum^2_{i=1} A_i x A^*_i.$ It is shown in
Section 5 of \cite{GKL06} (and not difficult to verify directly) that
this map is ergodic.
We will use the same notations here as in previous sections.
Observe that the vector
$\Omega_\CH:=\frac{1}{\sqrt{3}} (1,1,1)^T$
gives an invariant vector state for $Z$ as
$$\langle \Omega_\CH,Z(x)\Omega_\CH \rangle =
\langle \Omega_\CH ,x \Omega_\CH \rangle = \frac{1}{3}\sum^3_{i,j=1} x_{ij}.$$
$A^*_i \Omega_\CH = \frac{1}{\sqrt{2}}  \Omega_\CH$ and
hence $\uomega=\frac{1}{\sqrt{2}}(1,1).$ The orthogonal complement $\Ho$
of $\C \Omega_\CH$ in $\C^3$ and the orthogonal projection $Q$ onto $\Ho$
are given by
$$\Ho=\{ \left(\begin{array}{c} k_1 \\ k_2 \\ -(k_1 + k_2)\end{array} \right):
k_1,k_2 \in \C\}, \qquad
Q = \frac{1}{3} \left( \begin{array}{ccc}
2 & -1 & -1 \\ -1 & 2 & -1 \\ -1 & -1 & 2 \end{array} \right).
$$
From this we get for $\Ao_i = Q A_i Q = A_i Q$
$$\Ao_1 = \frac{1}{3 \sqrt{2}} \left( \begin{array}{ccc}
0 & 0 & 0 \\ 2 & -1 & -1 \\ -2 & 1 & 1 \end{array} \right),
~~\Ao_2 = \frac{1}{3 \sqrt{2}} \left( \begin{array}{ccc}
1 & 1 & -2 \\ -1 & -1 & 2 \\ 0 & 0 & 0 \end{array} \right).$$
We notice that the tuple $\uAo=(\Ao_1,\Ao_2)$ is $*$-stable as (by induction)
$$\sum_{|\alpha|=n} \Ao_\alpha \Aso_\alpha
= \frac{1}{3\times 2^{n-1}} \left( \begin{array}{ccc}
1 & -1 & 0 \\ -1 & 2 & -1 \\ 0 &
-1 & 1 \end{array} \right) \to 0 \quad (n \to \infty).
$$
Here $\CP=\C^2$ and $\Po:=\CP \ominus \C \Omega_\CP$
with $\OP = \frac{1}{\sqrt{2}} (1,1)^T$.
Easy calculation shows that $\hat{D}_*:\; \Ho \to \Po$ is given by
$$\left( \begin{array}{c} k_1 \\ k_2 \\ -(k_1 + k_2)\end{array} \right)
\mapsto \frac{1}{\sqrt{6}} (2k_1 + k_2)\left( \begin{array}{c} -1 \\ 1
\end{array} \right).$$
Moreover
$\Dso=\frac{1}{\sqrt{6}} \left( \begin{array}{ccc}
1 & 0 & -1 \\ 0 & 0 & 0 \\ -1 & 0 & 1 \end{array} \right).$
There exists an isometry $\gamma:\,\CDso \to \Po$
such that $\left( \begin{array}{c} 1 \\ 0 \\ -1 \end{array} \right)
\mapsto \left( \begin{array}{c} -1 \\ 1
\end{array} \right)$ and $\gamma(\Dso h) = \hat{D}_* h$ for $h \in \Ho$.

The map $\hat{C}: \CH \to \Gamma (\C^d) \otimes \CD_{\omega}$ is given by
$\hat{C}(\Omega_\CH)=1$ and for $h \in \Ho$ by
\begin{eqnarray*}
\hat{C} \left( \begin{array}{c} k_1 \\ k_2 \\ -(k_1 + k_2)\end{array} \right)
&=& e_0 \otimes \frac{(2k_1 + k_2)}{\sqrt{6}}
\left( \begin{array}{c} -1 \\ 1 \end{array}
\right) + \sum_{\alpha,\alpha_1=1} e_\alpha \otimes (\frac{1}{\sqrt{2}})^{|\alpha|} \frac{(k_1 + 2k_2)}{\sqrt{6}} \\
&& \times \left( \begin{array}{c} -1 \\ 1 \end{array} \right)
+ \sum_{\alpha,\alpha_1=2} e_\alpha \otimes
(\frac{1}{\sqrt{2}})^{|\alpha|} \frac{(k_1 - k_2)}{\sqrt{6}}
\left( \begin{array}{c} -1 \\ 1 \end{array} \right)
\end{eqnarray*}
where the summations are taken over all $0 \not=\alpha \in \tilde{\Lambda}$ such that $\alpha_i \neq \alpha_{i+1}$ for all $1 \leq i \leq |\alpha|$ and fixing $\alpha_1$ to $1$ or
$2$ as indicated. This simplification occurs because $\Ao^2_i = 0$ for
$i=1,2$.
All the summations below in this section are also of the same kind.

Now using the equations 4.2 and 4.3 for $\hat{\theta}_A :\CD_A \to \Gamma (\C^d) \otimes \CD_{\omega}$ and simplifying we get
\begin{eqnarray*}
\hat{\theta}_A\, d^1_{\Omega_\CH}  &=& -e_0 \otimes \frac{1}{6}
\left( \begin{array}{c} -1 \\ 1 \end{array}
\right) + \sum_{\alpha,\alpha_1=1} e_\alpha \otimes
(\frac{1}{\sqrt{2}})^{|\alpha|} \frac{1}{6}
\left( \begin{array}{c} -1 \\ 1 \end{array} \right)\\
&&+ \sum_{\alpha,\alpha_1=2} e_\alpha \otimes
(\frac{1}{\sqrt{2}})^{|\alpha|} \frac{1}{6}
\left( \begin{array}{c} -1 \\ 1 \end{array} \right)
= - \hat{\theta}_A\, d^2_{\Omega_\CH},
\end{eqnarray*}
and for $h \in \Ho,$
\begin{eqnarray*}
\hat{\theta}_A d^1_h
&=& -e_0 \otimes \frac{k_1}{2 \sqrt{3}} \left( \begin{array}{c} -1 \\
1 \end{array}
\right)+ e_1 \otimes \frac{(k_1 + k_2)}{\sqrt{6}}
\left( \begin{array}{c} -1 \\ 1 \end{array} \right)
+ \sum_{\alpha,\alpha_1=1} e_1 \otimes e_\alpha \\
&&\otimes (\frac{1}{\sqrt{2}})^{|\alpha|} \frac{(k_1 + 2k_2)}{\sqrt{6}}
\left( \begin{array}{c} -1 \\ 1 \end{array} \right)
- \sum_{\alpha,\alpha_1=2} e_1 \otimes e_\alpha \otimes
(\frac{1}{\sqrt{2}})^{|\alpha|} \frac{ k_2}{\sqrt{6}}
\left( \begin{array}{c} -1 \\ 1 \end{array} \right)\\
&&+ \sum_{\alpha,\alpha_1=2} e_\alpha \otimes
(\frac{1}{\sqrt{2}})^{|\alpha|} \frac{k_1}{2 \sqrt{3}}
\left( \begin{array}{c} -1 \\ 1 \end{array} \right),
\end{eqnarray*}

\begin{eqnarray*}
\hat{\theta}_A d^2_h&=& -e_0 \otimes \frac{(k_1 + k_2)}{2 \sqrt{3}}
\left( \begin{array}{c} -1 \\ 1 \end{array}
\right) + \sum_{\alpha,\alpha_1=1} e_\alpha \otimes
(\frac{1}{\sqrt{2}})^{|\alpha|} \frac{(k_1 + k_2)}{2 \sqrt{3}}
\left( \begin{array}{c} -1 \\ 1 \end{array}
\right)\\
&& + e_2 \otimes \frac{k_1}{\sqrt{6}} \left( \begin{array}{c} -1 \\ 1
\end{array} \right) + \sum_{\alpha,\alpha_1=1} e_2 \otimes e_\alpha
\otimes (\frac{1}{\sqrt{2}})^{|\alpha|}\frac{k_2}{\sqrt{6}}
\left( \begin{array}{c} -1 \\ 1 \end{array}
\right)\\
&& + \sum_{\alpha,\alpha_1=2} e_2 \otimes e_\alpha \otimes
(\frac{1}{\sqrt{2}})^{|\alpha|} \frac{(k_1 - k_2)}{\sqrt{6}}
\left( \begin{array}{c} -1 \\ 1 \end{array} \right).
\end{eqnarray*}

Form this we can easily obtain $\Co$ and $\theta_{\Ao}$ for $h \in \Ho$
by using the following relations from Theorem 5.2,
$$(\eins \otimes \gamma) \Co h = \hat{C} h,$$
$$(\eins \otimes \gamma) \theta_{\Ao} d^i_h = \hat{\theta}_A d^i_h.$$
Further
$$l_1= A_1 \Omega_\CH - \frac{1}{\sqrt{2}} \Omega_\CH=  \frac{1}{\sqrt{6}}
\left( \begin{array}{c} -1 \\ 0 \\ 1 \end{array} \right),~~
l_2= A_2 \Omega_\CH - \frac{1}{\sqrt{2}} \Omega_\CH=  \frac{1}{\sqrt{6}}
\left( \begin{array}{c} 1 \\ 0 \\ -1 \end{array} \right),$$
$\Ao_1 l_1=  \frac{1}{2 \sqrt{3}}  \left( \begin{array}{c} 0 \\ -1 \\ 1
\end{array} \right)$ and  clearly
$\Ho=\overline{\mbox{span}} \{ \Ao_\alpha l_i: i=1,2 \mbox{~and~}
\alpha \in \tilde{\Lambda}\},$ as already observed in Remark 3.1.



\begin{thebibliography}{1}

\bibitem [Ar03]{Ar03}
Arveson, W.:
{\em Noncommutative Dynamics and E-Semigroups,\/}
Springer Monographs in Mathematics (2003)

\bibitem [BP94]{BP94}
Bhat, B. V. R.; Parthasarathy, K. R.: \textit{Kolmogorov's existence
theorem for Markov processes in $C^*$ algebras}, Proc. Indian Acad.
Sci., Math. Sci., \textbf{104} (1994), 253-262

\bibitem [BJKW00]{BJKW00}
Bratteli, O.; Jorgensen, P.; Kishimoto, A.; Werner, R. F.:
\textit{Pure states on $\CO_d,$} J. Operator Theory, \textbf{43} (2000), 97-143

\bibitem [BJP96]{BJP96}
Bratteli, O.; Jorgensen, P. E. T.; Price, G. L.:
\textit{Endomorphisms of $B(\CH)$.}
Quantization, nonlinear partial differential equations, and operator algebra
(Cambridge, MA, 1994), 93-138,
Proc. Sympos. Pure Math.,\textbf{59}, Amer. Math. Soc., Providence, RI (1996).

\bibitem [Cu77]{Cu77}
Cuntz, J.: \textit{Simple $C^*$-algebras generated by isometries},
Commun. Math. Phys., \textbf{57} (1977),173-185

\bibitem [FF90]{FF90}
Foias, C.; Frazho, A. E.:
\textit{The commutant lifting approach to interpolation problems,}
Operator Theory: Advances and Applications, 44, Birkh\"auser
Verlag, Basel (1990).  

\bibitem [Go04]{Go04}
Gohm, R.: \textit{Noncommutative stationary processes,}
Lecture Notes in Mathematics, 1839, Springer-Verlag, Berlin (2004).

\bibitem [GKL06]{GKL06}
Gohm, R.; K\"ummerer, B.; Lang, T.:
\textit{Noncommutative symbolic coding},
Ergod.Th. \& Dynam.Sys., \textbf{26} (2006), 1521-1548

\bibitem [KM00]{KM00}
K\"ummerer, B.; Maassen, H.:
\textit{A scattering theory for Markov chains}.
Inf. Dim. Analysis, Quantum Prob. and Related Topics,
vol.\textbf{3} (2000), 161-176

\bibitem [K\"{u}85]{Ku85}
K\"{u}mmerer, B.: \textit{Markov dilations on $W^*$-algebras},
J. Funct. Anal., \textbf{63} (1985), 139-177

\bibitem [NF70]{NF70}
Sz.-Nagy, B., Foias, C.:
\textit{Harmonic analysis of operators on Hilbert space,}
North Holland Publ., Amsterdam-Budapest (1970)

\bibitem [Po89a]{Po89a}
Popescu, G.:
\textit{Isometric dilations for infinite sequences of noncommuting operators},
Trans. Amer. Math. Soc., \textbf{316} (1989), 523-536

\bibitem [Po89b]{Po89b}
Popescu, G.: \textit{Characteristic functions for infinite
sequences of noncommuting operators}, J. Operator Theory,
\textbf{22} (1989), 51-71

\bibitem [Po89c]{Po89c}
Popescu, G.: \textit{Multi-analytic operators and some
factorization theorems}, Indiana Univ.Math.J.,
\textbf{36} (1989), 693-710

\bibitem [Po95]{Po95}
Popescu, G.: \textit{Multi-analytic operators on Fock spaces}
J. Funct. Anal. \textbf{161} (1999), 27-61

\bibitem [Pow88]{Pow88}
Powers, R. T.:
\textit{An index theory for semigroups of $*$-endomorphisms of $B(\CH)$ and
type $II_1$ factors}, Canad. J. Math. \textbf{40} (1988),  no. 1,
86-114


\end{thebibliography}
\end{document}